\documentclass[11pt,letterpaper]{article}

\usepackage[utf8]{inputenc}
\usepackage[T1]{fontenc}
\usepackage[english]{babel}
\usepackage{lmodern}
\usepackage[margin=1in]{geometry}
\usepackage{microtype}

\usepackage{amsmath,amssymb,amsthm}
\usepackage{mathtools}
\usepackage{mathrsfs}          
\usepackage{bm}               

\numberwithin{equation}{section}

\usepackage{graphicx}
\usepackage{xcolor}
\usepackage{booktabs}
\usepackage{array}
\usepackage{enumitem}
\usepackage{tikz}
\usetikzlibrary{arrows.meta,positioning,calc}

\usepackage{titlesec}

\titleformat{\section}
  {\normalfont\Large\scshape}
  {\thesection}{1em}{}
\titleformat{\subsection}
  {\normalfont\large\scshape}
  {\thesubsection}{1em}{}
\titleformat{\subsubsection}
  {\normalfont\normalsize\scshape}
  {\thesubsubsection}{1em}{}

\titlespacing*{\section}
  {0pt}{3.5ex plus 1ex minus .2ex}{2.3ex plus .2ex}
\titlespacing*{\subsection}
  {0pt}{3.25ex plus 1ex minus .2ex}{1.5ex plus .2ex}
\titlespacing*{\subsubsection}
  {0pt}{3.25ex plus 1ex minus .2ex}{1.5ex plus .2ex}

\renewenvironment{abstract}
  {\par\vspace{1.5em}%
   \begin{center}\normalfont\large\abstractname\end{center}%
   \vspace{0.5em}\normalfont\normalsize\noindent\ignorespaces}
  {\par\vspace{1em}}

\usepackage{tocloft}
\newtheoremstyle{budapestplain}
  {\topsep}{\topsep}
  {\itshape}                  
  {}                          
  {\normalfont\scshape}       
  {.}{.5em}{}

\newtheoremstyle{budapestdefinition}
  {\topsep}{\topsep}
  {\normalfont}               
  {}
  {\normalfont\scshape}
  {.}{.5em}{}

\theoremstyle{budapestplain}

\newtheorem{theorem}{Theorem}[section]
\newtheorem{lemma}[theorem]{Lemma}
\newtheorem{proposition}[theorem]{Proposition}
\newtheorem{corollary}[theorem]{Corollary}

\newtheorem*{theorem*}{Theorem}
\newtheorem*{lemma*}{Lemma}

\theoremstyle{budapestdefinition}

\newtheorem{remark}[theorem]{Remark}

\newtheorem*{remark*}{Remark}
\newtheorem*{example*}{Example}
\newtheorem*{notation*}{Notation}

\usepackage[colorlinks=true,linkcolor=blue,citecolor=blue,urlcolor=blue]{hyperref}
\usepackage[nameinlink,noabbrev]{cleveref}
\crefname{exercise}{exercise}{exercises}
\Crefname{exercise}{Exercise}{Exercises}

\newcommand{\R}{\mathbb{R}}
\newcommand{\N}{\mathbb{N}}
\newcommand{\Z}{\mathbb{Z}}

\DeclarePairedDelimiterX{\inner}[2]{\langle}{\rangle}{#1,#2}

\title{Exceptional parameters in the Erd\H{o}s similarity problem\\[0.2em]}
\author{N. Kulkarni and N. Mora Cu\'ellar}
\date{\today} 

\begin{document}

\maketitle

\begin{abstract}
The Erd\H os similarity conjecture says that for every infinite set $A \subset \mathbb{R}$, there exists a set $E$ of positive Lebesgue measure which does not contain any affine copy of $x + tA$ of $A$ (with $t \neq 0$). Kolountzakis proved that for every infinite $A$, there exists $E$ of positive Lebesgue measure so that the exceptional sets $D_E(A)$ of
scales $t$, and $\mathcal E_E(A)$ of affine maps $(x,t)$, for which
$x+tA\subseteq E$, are both Lebesgue null. We show that these exceptional
parameter sets are small in two finer senses. First, given any countable
family of Fourier decay rates, a single $E$ can be constructed so that $D_E(A)$ is null with respect to every
measure whose Fourier transform decays at one of these rates. Consequently, there are compact sets $K$ of measure arbitrarily close to $1$
for which $D_K(A)$ and $\mathcal E_K(A)$ have modified Fourier dimension
zero. Second, we give a criterion,
depending both on the growth of a gauge function $h$ and on the gaps between elements of lacunary sequences in $A$, which guarantees that $\mathcal H^h(D_E(A))=0$ for gauges
with $h(u)/u\to\infty$, as well as a similar result for
$\mathcal E_K(A)$. For sets such as $A = \{1/n!\}$ and $A = \{e^{-n^q}\}$,
this yields explicit gauges of the form $h(u) = u\exp((\log 1/u)^\theta)$ with $0 < \theta < 1$.
\end{abstract}

\section{Introduction and results}
\addcontentsline{toc}{section}{Introduction}
The Erd\H{o}s similarity conjecture \cite{Erdos,JLM} asserts that for
every infinite set $A\subseteq\mathbb R$ there is a set $E$ of positive
Lebesgue measure that contains no affine copy $x+tA$ with $t\neq0$. 
The
conjecture remains open in general, but Kolountzakis
\cite[Theorem~2]{Kolountzakis} proved a measure theoretic version of it. He showed that for every
infinite $A$, there exists a set $E$ of positive Lebesgue measure such that, for almost
every scale $t$, no translate of $tA$ lies in $E$. 
Equivalently, the
affine copies of $A$ inside $E$ form a Lebesgue-null family of
parameters in $\mathbb{R}^2$. The goal of this paper is to give quantitative refinements of Kolountzakis' result.
To state these refinements, we first define, for $A\subseteq\mathbb R$ and a closed set
$E\subseteq\mathbb R$, the set of scales and the set of affine maps that place
a copy of $A$ inside $E$:
\[
 \begin{aligned}
 D_E(A)&=\{t\ne0:x+tA\subseteq E\text{ for some }x\in\mathbb R\},\\
 \mathcal E_E(A)&=\{(x,t):t\ne0,\ x+tA\subseteq E\}.
 \end{aligned}
\]
Kolountzakis's theorem
constructs a set $E$ of positive measure with $\mathcal L^1(D_E(A))=0$ and
$\mathcal L^2(\mathcal E_E(A))=0$, where $\mathcal L^d$ denotes Lebesgue measure on $\mathbb R^d$. We will show that (1)
these exceptional parameter sets 
are null with respect to any measure whose Fourier transform decays at a prescribed
rate, and (2) they are null with respect to Hausdorff measures associated to gauges $h(u)$ which we will specify precisely. In both cases, a single set $E$ works simultaneously
for a countable family of measures.

The Fourier
transform of a finite Borel measure $\nu$ on $\mathbb R^d$ is the integral
$\widehat\nu(\xi)=\int e^{-2\pi i\xi\cdot y}\,d\nu(y)$.
A \emph{decay modulus} is a non-increasing function
$\omega:[1,\infty)\to(0,1]$ with $\omega(r)\to0$ as $r\to\infty$. A finite
Borel measure $\nu$ on $\mathbb R$ has \emph{decay $\omega$} if 
\[
 |\widehat\nu(\xi)|\lesssim \omega(|\xi|)\qquad\text{for all }|\xi|\ge1.
\]
Our first result makes $D_E(A)$ null for all measures with a given
countable family of decay moduli at once.

\begin{theorem}\label{thm:prescribed-moduli}
Let $A\subset \mathbb R$ be an infinite set, and let $(\omega_m)_{m\ge1}$ be a
countable family of decay moduli. For every $\varepsilon>0$ there is a
closed, $1$-periodic, nowhere-dense set $E\subseteq\mathbb R$ such that
$\mathcal L^1(E\cap I)>1-\varepsilon$ for every interval $I$ of length $1$, and
\[
 \nu(D_E(A))=0
\]
for every finite positive Borel measure $\nu$ with decay $\omega_m$ for
some $m\geq 1$.
\end{theorem}

Since each decay modulus controls the uncountable class of all measures with $|\widehat\nu|\lesssim \omega$, Theorem $\ref{thm:prescribed-moduli}$ allows us to control the modified Fourier dimension of the exceptional sets.
For a finite Borel measure $\nu$ on $\mathbb R^d$, its Fourier dimension is defined by
\[
 \dim_F\nu=\sup\{s\in[0,d]:
 |\widehat\nu(\xi)|\lesssim (1+|\xi|)^{-s/2}\text{ for every }\xi\}.
\]
For a Borel
set $B\subseteq\mathbb R^d$, its modified Fourier dimension
\cite{EkstromPerssonSchmeling} is then defined by first setting $\dim_{FM}\varnothing=0$ and for nonempty sets:
\[
  \dim_{FM}B=\sup\{\dim_F\nu:\nu\text{ is a Borel probability measure and } \nu(B)>0\}.
\]
If $\dim_{FM}B=0$, this tells us that $B$ is null with respect to every probability measure with some positive power of Fourier decay.
\begin{theorem}
\label{thm:exceptional-main}
For every infinite $A \subset \mathbb R$ and every $\varepsilon>0$, there is a compact
nowhere-dense $K\subseteq[0,1]$ with $\mathcal L^1(K)>1-\varepsilon$ such that
\[
 \mathcal L^1(D_K(A))=0,\qquad
 \dim_{FM}D_K(A)=0,\qquad
 \dim_{FM}\mathcal E_K(A)=0.
\]
Moreover, $\nu(D_K(A))=0$ for every finite positive Borel measure $\nu$
with $|\widehat\nu(\xi)|\lesssim (\log(2+|\xi|))^{-p}$ at all
sufficiently large frequencies, for some $p>0$.
\end{theorem}


We next turn to Hausdorff measures. 
The conclusion of Kolountzakis' theorem, i.e. $\mathcal L^1(D_K(A))=0$, is equivalent to the statement
$\mathcal H^1(D_K(A))=0$. Below we give a refinement of this statement using Hausdorff measures associated to gauge functions. For any bounded and infinite set $A$, we will describe how to construct sets $E$ and $K$ of positive Lebesgue measure and gauge functions $h$ and $g$ such that $\mathcal{H}^h(D_E(A)) = 0$ and $\mathcal{H}^g(\mathcal{E}_K(A)) = 0$. The chosen gauge function $h$ will also satisfy
\[\lim_{u \to 0} \frac{h(u)}{u} = \infty,\]
so $\mathcal H^h(F)=0$ implies $\mathcal H^1(F)=0$, while there are sets with
$\mathcal H^1(F)=0<\mathcal H^h(F)$. Thus, $\mathcal H^h(D_E(A))=0$ is a strict
refinement of Kolountzakis' result. In addition to proving the existence of such a gauge function $h$, our proof gives a precise criterion based on the structure of $A$ for gauge functions $h$ to satisfy $\mathcal{H}^h(D_E(A)) = 0$.

We begin by defining a gauge function, its associated Hausdorff outer measure, and introduce the concept of the \emph{overhead} of a gauge function.
A \emph{gauge} is a nondecreasing function $h:(0,1]\to(0,\infty)$ with
$h(u)\to0$ as $u\to 0$. The associated \emph{Hausdorff outer measure}
is
\[
 \mathcal H^h(F)=\lim_{\delta\downarrow0}
 \inf\big\{\sum_i h(\operatorname{diam}U_i):
 F\subseteq\bigcup_iU_i,\quad 0<\operatorname{diam}U_i\le\delta\big\},
\]
where the infimum runs over countable covers $\{U_i\}$ of $F$. The \emph{overhead} of $h$ is
the nondecreasing function
\begin{equation}
 \Psi_h(S)=
 \sup_{0\le s\le S}
 \big( \log \left(e^{s}h(e^{-s})\right)\big)_+,
 \qquad S\ge0.
 \label{SN:eq:Psi-definition}
\end{equation}
By construction,
\begin{equation}\label{SN:eq:overhead-bound}
 h(u)\le u\,e^{\Psi_h(\log(1/u))}.
\end{equation}
The relevant gauges are the ``almost linear'' ones, where $\Psi_h(S)=o(S)$ and $h(u)\le u^{1-o(1)}$.
Two model examples are
\begin{align}
 h_\theta(u)&=u\exp\!\left((\log(1/u))^\theta\right),
       &&0<\theta<1, \label{SN:eq:h-theta}\\
 \ell_p(u)&=u(\log(1/u))^p,
       &&p>0, \label{SN:eq:ell-p}
\end{align}
and their overheads are $\Psi_{h_\theta}(S)=S^\theta+O(1)$ and
$\Psi_{\ell_p}(S)=p\log S+O(1)$.

The Erd\H os similarity conjecture has already been solved for geometric sequences \cite{Geometric}, and Jiang and Xi \cite{JiangXi} extended this result to every positive zero sequence $(b_n)$ with bounded logarithmic gaps $g_n:=\log(b_n/b_{n+1})$.
Their results do not apply to patterns for which $g_n$ is unbounded, such as the sequences $\{1/n!\}$ or $\{e^{-n^2}\}$.
The following theorem applies to every bounded infinite $A$. The allowed gauges $h$ depend on the growth, as $n \to \infty$, of  $\log(b_1/b_n)=g_1+\cdots+g_{n-1}$, where $b_1>\cdots>b_n$ is a lacunary block contained in $A$. The faster the gaps grow, the closer the gauge must be to a linear function, and the weaker the conclusion of the theorem.

\begin{theorem}
\label{SN:thm:block-entropy}
Let $A\subseteq\mathbb R$ be bounded and infinite, and let
$(h_m)_{m\ge1}$ be a countable family of gauges. Suppose that for every
$m\ge1$ and every $\rho\in(0,1)$ there are arbitrarily large $n$, together with
$x\in\mathbb R$, $\sigma\in\{-1,1\}$ and $b_1>\cdots>b_n>0$, such that
\begin{equation}\label{SN:eq:block-criterion}
 x+\sigma\{b_1,\ldots,b_n\}\subseteq A,\qquad
 b_{i+1}\le\rho\,b_i\quad(1\le i<n),\qquad
 \Psi_{h_m}\Bigl(\frac1\rho\log\frac{b_1}{b_n}\Bigr)\le\rho\, n.
\end{equation}
Then for every $\varepsilon>0$ there is a closed $1$-periodic
nowhere-dense $E\subseteq\mathbb R$ such that
\[
 \mathcal L^1(E\cap I)>1-\varepsilon\quad\text{for every }|I|=1,\qquad
 \mathcal H^{h_m}(D_E(A))=0\quad\text{for every }m\ge1.
\]
\end{theorem}
\noindent Any gauge bound on $D_K(A)$ implies a gauge bound on the set of exceptional affine maps $\mathcal{E}_K(A)$ through a companion gauge.

\begin{proposition}\label{SN:prop:pairs}
Let $A$ contain at least two points, let $K\subseteq[0,1]$ be compact,
and let $h$ be a Hausdorff gauge. If $\mathcal{H}^h(D_K(A))=0$, then
\[
 \mathcal{H}^{g_h}(\mathcal{E}_K(A))=0,\qquad \text{ where }\qquad g_h(u)=uh(u/2).
\]
In particular, \eqref{SN:eq:h-theta} and \eqref{SN:eq:ell-p} give
\[
g_\theta(u)= u^2\exp\!\left((\log(1/u))^\theta\right),\qquad
 g_p(u) = u^2(\log(1/u))^p.
\]
\end{proposition}

\noindent The paper is organized as follows. Section~2 proves
Theorems~\ref{thm:prescribed-moduli} and~\ref{thm:exceptional-main}. Section~3 proves
Theorem~\ref{SN:thm:block-entropy} and Proposition~\ref{SN:prop:pairs},
and gives explicit examples of patterns and gauges to which
Theorem~\ref{SN:thm:block-entropy} applies.

\section{Almost no affine copy with respect to measures with Fourier decay}

In this section we prove Theorems~\ref{thm:prescribed-moduli}
and~\ref{thm:exceptional-main}.
Given a bounded and infinite set $A\subset \R$ and $\varepsilon>0$, we will construct the set $E$ by simultaneously constructing  sequences of Borel sets $\{H_j\}_j$ and $\{B_j\}_j$ such that
\begin{itemize}
    \item $\sum_j\mathcal L^1(H_j\cap I)=\varepsilon$ for any interval $I$ of length 1;
    \item $D_{H_j^c}(A)\subset B_j$ for every $j\geq 1$; and
    \item $\nu(B_j)\leq 2^{-j}$ for any finite measure $\nu$ with a decay modulus in the family $\{\omega_m\}_m$.
\end{itemize}
Finally, we will take $E=\bigcap_j H_j^c$, so $\mathcal L^1(E\cap I)\geq 1-\varepsilon$ while $\nu(D_E(A))=0$.
To construct these sequences of sets we need some quantitative lemmas, proved in the following subsection.
\subsection{Finite blocks and periodic holes}
We use the following lemma for both Fourier testing and
Hausdorff-content estimates. 
Write $\mathbb T=\mathbb R/\mathbb Z$ and $\|V\|_{A(\mathbb T)}=\sum_k|\widehat V(k)|$ for the Wiener norm.
We will construct a family of trigonometric polynomials that are majorants for a family of indicator functions $\{\chi_{I_\ell^c}\}_{\ell=1}^N$, with integral uniformly close to $1$, and with degrees and norms uniformly bounded. The proof uses standard facts about the Fej\'er kernel (nonnegativity, normalization, degree, and the tail estimate) which can be found in books such as \cite[Section~4.1]{Montgomery2014}.

\begin{lemma}[Majorants of fixed degree]
\label{SN:lem:atlas}
For every $0<\alpha<1/8$ there are closed arcs
$I_1,\ldots,I_N\subseteq\mathbb T$ and nonnegative trigonometric polynomials
$V_1,\ldots,V_N$ such that:
\begin{itemize}
\item every open arc of length $\alpha$ contains some $I_\ell$;
\item $N\le C/\alpha$;
\item $V_\ell\ge1$ on $\mathbb T\setminus I_\ell$;
\item $\int_\mathbb T V_\ell\le1-c\alpha$; and
\item all degrees are at most a finite $R=R(\alpha)$ and all Wiener
      norms are at most a finite $W=W(\alpha)$.
\end{itemize}
\end{lemma}

\begin{proof}
Choose $N=\lceil4/\alpha\rceil$ equally spaced centers and let
$I_\ell$ be the closed arc of length $\alpha/2$ about the $\ell$th
center. The space between consecutive centers is at most $\alpha/4$, so every open arc of
length $\alpha$ contains one.  Inside $I_\ell$,
let $J_\ell$ be the central arc of length $\mu:=\alpha/8$.  If $F_R$
is the normalized Fej\'er kernel of order $R$, i.e.
\[F_R(y) = \sum_{|k|\leq R} \big(1-\tfrac{|k|}{R+1}\big) e^{2\pi i ky},\]
set $U_\ell=1_{J_\ell}*F_R$.
Since $F_R$ is an approximation to the identity, we have $0\le U_\ell\le1$ and
$\int_\mathbb T U_\ell=|J_\ell|=\mu$.
Using the standard Fejer tail estimate (\cite[Section ~4.1]{Montgomery2014}), noting that $\operatorname{dist}(J_\ell, I_\ell^c)\leq 3\alpha/16$, and taking $R=R(\alpha)$ to be sufficiently large, we can guarantee that
\[
 \sup_{\mathbb T\setminus I_\ell}U_\ell\le \frac{1}{(R+1)\sin^2(3\pi \alpha/16)} \leq \mu/8.
\]
Then define $V_\ell :=\tfrac{1-U_\ell}{1-\mu/8}$,
so that $V_\ell \geq 1$ outside $I_\ell$ and
\[
 \int_\mathbb T V_\ell
 =\frac{1-\mu}{1-\mu/8}
 \le1-7\alpha/64.
\]
Since each $V_\ell$ has only finitely many nonzero Fourier coefficients, it has a finite Wiener norm. There are only $N$ trigonometric polynomials $V_\ell$, so there is a uniform bound $W(\alpha)$ for all $||V_\ell||_{A(\mathbb{T})}$.
\end{proof}

We may assume that the set $A$ is bounded (otherwise, the Erd\H{o}s similarity conjecture holds trivially), so $A$ has an accumulation point, say $a_*$. Thus, we can always find a positive sequence $b_i\downarrow 0$ such that $$a_*+\sigma b_i\in A$$ for $\sigma\in\{-1,1\}$ and for all $i\geq 1$. Passing to a subsequence, we may assume that it is as lacunary as we like. The idea is to bound the measures $\nu(B_j)$ in terms of the Fourier structure of such lacunary blocks in $A$. 

\begin{lemma}\label{lemma:tail-domination}
Let $R\ge1$ be an integer and let $b_1>\cdots>b_n>0$ satisfy
$b_{i+1}\le b_i/(4R)$ for $1\le i<n$. Then $\sum_{i=1}^nb_i\le\frac43b_1$ and
\[
 \Bigl|\sum_{i=1}^nk_ib_i\Bigr|\ge b_n
 \qquad\text{for every }k\in\mathbb Z^n\setminus\{0\}\text{ with }|k_i|\le R.
\]
\end{lemma}

\begin{proof}
The sum $\sum_{i=1}^nb_i$ can be majorized by a geometric series with first term $b_1$ and common ratio $1/4$, so this geometric series equals $4b_1/3$.
Next, let $k\in\mathbb{Z}^n\setminus\{0\}$ with $|k_i|\le R$, and let $i_0=\min\{i:k_i\ne0\}$. If
$i_0=n$, then $|\sum_ik_ib_i|=|k_n|b_n\ge b_n$. Otherwise,
\[
 \Bigl|\sum_{i=1}^n k_ib_i\Bigr|
 \ge b_{i_0}-R\sum_{i>i_0}b_i
 \ge b_{i_0}\Bigl(1-R\cdot\tfrac{(1/4R)}{1-(1/4R)}\Bigr)
 =b_{i_0}\Bigl(1-\tfrac{R}{4R-1}\Bigr)
 \ge\tfrac23b_{i_0}
 \ge\tfrac23\cdot4R\,b_n\ge b_n,
\]
since $b_{i_0}\ge b_{n-1}\ge4Rb_n$.
\end{proof}
The following lemma constructs sets that contain some dilation parameters and whose measure $\nu$ is bounded in terms of its Fourier transform $\hat{\nu}$.
\begin{lemma}\label{lem:missing-arc}
Fix $0<\alpha<1/8$, and let $R=R(\alpha)$ be as in Lemma~\ref{SN:lem:atlas}. There are constants $k_\alpha>0$, $C_\alpha<\infty$ and $W_\alpha\ge1$, depending only on $\alpha$, with the following property. Let $q>0$, and let $P=(b_1>\cdots >b_n>0)$ satisfy $b_{i+1}\le b_i/(4R)$ for $1\le i<n$.
Then there is an open set $B(P,q)\subseteq\mathbb R$ such that:
\begin{itemize}
    \item if the chain $\{tqb_i \mod 1, 1\leq i\leq n\}$ has a missing open arc of length $\alpha$, then $t\in B(P,q)$;
    \item for every finite positive Borel measure $\nu$,
    \begin{equation}
         \nu(B(P,q))\leq C_{\alpha}\left[\nu(\R)e^{-k_{\alpha}n}+W_\alpha^n\sup_{|\xi|\geq qb_n}|\hat{\nu}(\xi)|\right];
         \label{eq:measure-bound}
    \end{equation}

    \item For every $T>0$, $B(P,q)\cap[-T,T]$ has at most
$$\label{eq:common-component-bound}
 C_\alpha\left(n+Tq\sum_i b_i\right)
$$
interval components.
\end{itemize}
\end{lemma}

   \begin{proof}
Let $I_1,\ldots,I_N$ and $V_1,\ldots,V_N$ be as in Lemma~\ref{SN:lem:atlas}, and set $W_\alpha=\max\{1,W(\alpha)\}$, with $W(\alpha)$ also from Lemma~\ref{SN:lem:atlas}. For each $\ell$, set
\[
 B_\ell=\{t:qb_i t\bmod1\notin I_\ell\ \text{for all }i\},
 \qquad
 B(P,q)=\bigcup_{\ell=1}^N B_\ell.
\]
Every missing arc of length $\alpha$ contains some $I_\ell$, so the
first property holds. Observe that
\[1_{B_\ell}(t)\le\prod_i V_\ell(qb_it).\]
Since $V_\ell$ is a trigonometric polynomial of degree $R$, expanding
the product gives
\[
 \prod_{i=1}^n V_\ell(qb_it)
 =\sum_{\substack{k\in\mathbb{Z}^n\\|k_i|\le R}}
 \Bigl(\prod_i\widehat V_\ell(k_i)\Bigr)
 e^{2\pi iqt\bigl(\sum_i k_ib_i\bigr)}.
\]
Integrating against $\nu$, the $k=0$ term contributes
$\widehat V_\ell(0)^n\nu(\mathbb{R})\le(1-c\alpha)^n\nu(\mathbb{R})$, from Lemma~\ref{SN:lem:atlas}.
For $k\ne0$, Lemma~\ref{lemma:tail-domination} gives
$|\sum_i k_ib_i|\ge b_n$, so every such frequency is
bounded by $\sup_{|\xi|\ge qb_n}|\widehat\nu(\xi)|$,
and the total weight of these terms is at most $W_\alpha^n$ since the sum
of absolute Fourier coefficients factorises as
$\prod_i\sum_{|k_i|\le R}|\widehat V_\ell(k_i)|\le W_\alpha^n$.
Setting $k_\alpha=-\log(1-c\alpha)$ and summing over the $N=N(\alpha)$
arcs proves the second property. For the third, the boundary
of the condition $qb_it\bmod1\notin I_\ell$ in $[-T,T]$ consists of
$O(1+Tqb_i)$ points, and summing over all $i$ and $\ell$ gives the
required bound.
\end{proof}

\subsection{Proof of the theorems}
\begin{proof}[Proof of Theorem~\ref{thm:prescribed-moduli}]
Now we put all of these parameters and constructions together.
Let $\varepsilon>0$. As we mentioned before, since $A$ is bounded an infinite, there is $a_*\in A$ accumulation point, and a sequence $\{b_i\}_i\downarrow 0$ such that $$a_*+\sigma b_i\in A$$ approach $a_*$ from one side. Here, $\sigma\in\{-1,1\}$ and without loss of generality we can assume that $\sigma=1$. Now choose a sequence $\sum_j\alpha_j=\varepsilon$, and for each $\alpha_j$ take $R_j, W_j, k_j$ and $C_j$ as in Lemma \ref{SN:lem:atlas} and Lemma \ref{lem:missing-arc}. Select $n_j\in\N$ so large that
$$C_je^{-k_jn_j}\leq 2^{-j},$$
and select $n_j$ terms $b_{j,1}>\cdots >b_{j,n_j}$ of the sequence $(b_i)$ with successive ratios at most $\frac{1}{4R_j}$; write $P_j=(b_{j,1},\ldots,b_{j,n_j})$. Finally, choose an integer $q_j\geq j$ such that $q_jb_{j,n_j}\geq j$ and such that for every $m\leq j$,

\begin{equation}
    C_jW_j^{n_j}\omega_m\left(q_jb_{j,n_j}\right)\leq 2^{-j}.
    \label{eq:j-bound}
\end{equation}

Now define $$H_j:=\left\{y:q_jy\mod 1\in\left(-\frac{\alpha_j}{2},\frac{\alpha_j}{2}\right)\right\}.$$
Each $H_j$ is $1-$periodic and has measure $\alpha_j$ inside each interval of length 1. So we can take
$$E=\bigcap_{j}H_j^c.$$
Notice that  $q_j^{-1}\Z\subset H_j$ for each $j\geq 1$, so $\bigcup_j H_j$ is dense and in consequence $E$ is nowhere dense, $1-$periodic and $\mathcal L^1(E\cap I)\geq 1-\sum_j\alpha_j=1-\varepsilon$ for any interval $I$ of length 1.

Also, if $t\in D_E(A)$, then there exists some $x\in \R$ such that $x+tA\subset E$. In particular, for any $j\geq 1$, and any $i\in\{1,\dots, n_j\}$, $$x+t(a_*+ b_{j,i})\subset H_j^c.$$
That is, $$q_jtb_{j,i}\mod 1\notin\left(\frac{-\alpha_j}{2}-q_jx,\frac{\alpha_j}{2}-q_jx\right).$$
So $\{q_jtb_{j,i}\mod 1: 1\leq i\leq n_j\}$ has a missing arc of length $\alpha_j$ and thus, by the construction made in Lemma \ref{lem:missing-arc}, $t\in B(P_j,q_j)$. In consequence, $D_E(A)\subset B(P_j,q_j)$ for any $j\geq 0$. Thus, putting together equations \ref{eq:measure-bound} and \ref{eq:j-bound}, we get that for any finite Borel measure $\nu$ with Fourier decay $\omega_m$,
$$\nu(D_E(A))\leq C_jW_j^{n_j}\omega_m\left(q_jb_{j,n_j}\right)\leq 2^{-j}\qquad \text{for any } j\geq m, $$
so $\nu(D_E(A)=0$.
\end{proof}
Theorem~\ref{thm:exceptional-main} follows from Theorem~\ref{thm:prescribed-moduli}.

   \begin{proof}[Proof of Theorem~\ref{thm:exceptional-main}]
We reduce to Theorem~\ref{thm:prescribed-moduli} by choosing a family of
decay moduli that simultaneously captures all positive power-decay rates,
all positive logarithmic-decay rates, and also the Fourier decay of Lebesgue
measure on bounded intervals. Specifically, define
\[
 \omega_{2m-1}(\xi)=(1+|\xi|)^{-1/m}
 \qquad\text{and}\qquad
 \omega_{2m}(\xi)=(\log(2+|\xi|))^{-1/m},\qquad m\ge1.
\]
Apply Theorem~\ref{thm:prescribed-moduli} to $A$, $\varepsilon$, and the family
$(\omega_m)_{m\ge1}$ to obtain a closed, $1$-periodic, nowhere-dense set
$E\subseteq\mathbb R$ with $\mathcal L^1(E\cap I)>1-\varepsilon$ for every unit interval $I$,
and $\nu(D_E(A))=0$ for every finite positive Borel measure $\nu$ having
decay $\omega_m$ for some $m\ge1$. Set $K=E\cap[0,1]$.

Since $E$ is closed and $1$-periodic, $K$ is compact and nowhere dense,
and $\mathcal L^1(K)>1-\varepsilon$. Any affine copy of $A$ contained in $K$ is also
contained in $E$, so $D_K(A)\subseteq D_E(A)$ and
$\mathcal{E}_K(A)\subseteq\mathcal{E}_E(A)$.

For $\mathcal L^1(D_K(A))=0$, fix any $T>0$ and let $\nu_T=\mathcal L^1|_{[-T,T]}$ be Lebesgue
measure restricted to $[-T,T]$.  As is well-known, its Fourier transform satisfies
\[
 |\hat\nu_T(\xi)|
 =\frac{|\sin(2\pi T\xi)|}{\pi|\xi|}
 \le\frac{1}{\pi|\xi|}
 \le\frac{2}{\pi}(1+|\xi|)^{-1}
 =\frac{2}{\pi}\,\omega_1(|\xi|)
\]
for all $|\xi|\ge1$, so $\nu_T$ has decay $\omega_1$. Theorem~\ref{thm:prescribed-moduli}
then gives $\nu_T(D_E(A))=0$, that is, $\mathcal L^1(D_E(A)\cap[-T,T])=0$.
Since $T$ was arbitrary, $\mathcal L^1(D_K(A))\le \mathcal L^1(D_E(A))=0$.

In order to prove that $\dim_{FM}D_K(A)=0$, let $\nu$ be any Borel probability measure
with $|\hat\nu(\xi)|\le C(1+|\xi|)^{-s/2}$ for some $s>0$. Choose any integer
$m\ge2/s$; then $1/m\le s/2$, and in consequence
\[
 |\hat\nu(\xi)|\le C(1+|\xi|)^{-s/2}\le C(1+|\xi|)^{-1/m}=C\,\omega_{2m-1}(|\xi|).
\]
So $\nu$ has decay $\omega_{2m-1}$, and Theorem~\ref{thm:prescribed-moduli}
gives $\nu(D_K(A))\le\nu(D_E(A))=0$. Since $\nu$ was an arbitrary Borel
probability measure with positive power Fourier decay, $\dim_{FM}D_K(A)=0$.

For $\dim_{FM}\mathcal{E}_K(A)=0$, let $\mu$ be a Borel probability measure
on $\mathbb R^2$ with $|\hat\mu(\bm\xi)|\le C|\bm\xi|^{-s/2}$ for some $s>0$.
Let $\nu=\pi_*\mu$ be its pushforward under the coordinate projection
$\pi(x,t)=t$. Then $\nu(D_K(A))\ge\mu(\mathcal{E}_K(A))$, and
\[
 |\hat\nu(\eta)|
 =\left|\int e^{-2\pi i\eta t}\,d\mu(x,t)\right|
 =|\hat\mu(0,\eta)|
 \le C|\eta|^{-s/2}.
\]
So $\nu$ is a finite positive measure with positive power Fourier decay, and
by the previous step $\nu(D_K(A))=0$. Thus $\mu(\mathcal{E}_K(A))=0$, and
since $\mu$ was arbitrary, $\dim_{FM}\mathcal{E}_K(A)=0$.

Finally, let $\nu$ be a finite positive Borel measure satisfying
$|\hat\nu(\xi)|\le C_\nu(\log(2+|\xi|))^{-p}$ for some $p>0$ and all
sufficiently large $|\xi|$. Choose an integer $m\ge1/p$; then $1/m\le p$,
and so
\[
 |\hat\nu(\xi)|
 \le C_\nu(\log(2+|\xi|))^{-p}
 \le C_\nu(\log(2+|\xi|))^{-1/m}
 =C_\nu\,\omega_{2m}(|\xi|)
\]
for all $|\xi|$ large enough. Since $\omega_{2m}$ is bounded below on
compact sets and $|\hat\nu|$ is bounded everywhere, increasing $C_\nu$
if necessary ensures the inequality holds for all $\xi$, so $\nu$ has decay
$\omega_{2m}$ and $\nu(D_K(A))\le\nu(D_E(A))=0$.
\end{proof}

\section{Almost no affine copy with respect to Hausdorff-gauge measures}
\subsection{Proof of the theorem}

It suffices to treat a single gauge $h$ satisfying
\eqref{SN:eq:block-criterion}: if $E_m$ is the set obtained for $h_m$ with
$\varepsilon2^{-m-1}$ in place of $\varepsilon$, then $E=\bigcap_mE_m$ is
closed, $1$-periodic and nowhere dense, satisfies
$\mathcal L^1(E\cap I)\ge1-\varepsilon/2$ for $|I|=1$, and
$D_E(A)\subseteq D_{E_m}(A)$ for every $m$.

\begin{proof}[Proof of Theorem~\ref{SN:thm:block-entropy}]
Fix a gauge $h$ satisfying \eqref{SN:eq:block-criterion}, and choose
$0<\alpha_j<1/8$ with $\sum_j\alpha_j<\varepsilon$. For each $j$ let
$R_j$, $W_j\ge1$, $k_j$ and $C_j\ge1$ be the constants of
Lemma~\ref{SN:lem:atlas} and Lemma~\ref{lem:missing-arc} for
$\alpha=\alpha_j$, and define
\[
 D_j:=\log W_j+2k_j+2+\log\bigl(C_j(1+2j)\bigr).
\]

\noindent \textbf{Step 1: Choose a block and a frequency.} Let
$\rho_j=\min\{1/(4R_j),\,(1+D_j/\log2)^{-1},\,k_j/2\}$. By
\eqref{SN:eq:block-criterion} with $\rho=\rho_j$, there are arbitrarily large
$n$ for which $A$ contains a set $x_j+\sigma_j\{b_{j,1},\ldots,b_{j,n}\}$
with $b_{j,i+1}\le\rho_j\,b_{j,i}$ and
$\Psi_h\bigl(\rho_j^{-1}\log(b_{j,1}/b_{j,n})\bigr)\le\rho_jn$. Fix such a block,
with length $n_j$ large enough that $2C_j(2j+1)e^{-k_jn_j/2}<2^{-j}$,
$C_j(2j+1)e^{-k_jn_j}\le2^{-j}$ and $e^{k_jn_j}\ge\operatorname{diam}A$, and write
\[
 P_j=(b_{j,1},\ldots,b_{j,n_j}),\qquad
 S_j=\log\frac{b_{j,1}}{b_{j,n_j}}.
\]
Since $\rho_j\le1/(4R_j)$, we have $S_j\ge(n_j-1)\log4\ge n_j\log2$,
so by our choice of $\rho_j$, 
\[S_j+D_jn_j\le \Bigl(1+\tfrac{D_j}{\log2}\Bigr) S_j\le \tfrac{S_j}{\rho_j}.\]
As $\rho_j\le k_j/2$,
the monotonicity of $\Psi_h$ and the theorem hypothesis gives
\begin{equation}\label{SN:eq:block-choice}
 \Psi_h\bigl(S_j+D_jn_j\bigr)\le\Psi_h\bigl(S_j/\rho_j\bigr)
 \le\rho_jn_j\le\frac{k_j}{2}n_j.
\end{equation}
Choose the frequency $q_j :=\bigl\lceil W_j^{n_j}e^{k_jn_j}/b_{j,n_j}\bigr\rceil$. Since
$b_{j,n_j}<b_{j,1}-b_{j,2}\le\operatorname{diam}A$ and $W_j\ge1$,
\begin{equation}\label{SN:eq:q-lower}
 q_j\ge\frac{W_j^{n_j}e^{k_jn_j}}{b_{j,n_j}}\ge\frac{e^{k_jn_j}}{\operatorname{diam}A}\ge\max\Bigl\{1,\frac{2^j}{\operatorname{diam}A}\Bigr\}.
\end{equation}
In particular $W_j^{n_j}e^{k_jn_j}/b_{j,n_j}\ge1$, so $q_j\le2W_j^{n_j}e^{k_jn_j}/b_{j,n_j}$. \\

\noindent\textbf{Step 2: Construct the set $E$.} Set
\[
 H_j=\bigl\{y:q_jy\bmod1\in(0,\alpha_j)\bigr\},
 \qquad
 E=\mathbb{R}\setminus\bigcup_j H_j.
\]
Explicitly, $H_j = \bigcup_{m\in\mathbb{Z}}\big(m/q_j, (m+\alpha_j)/q_j\big)$ so we observe that $H_j$ is open, $1$-periodic, and has $\mathcal L^1(H_j\cap I)=\alpha_j$ for every unit interval $I$. Thus, $E$ is closed and $1$-periodic with
$\mathcal L^1(E\cap I)\ge1-\sum_j\alpha_j>1-\varepsilon$. Every point of
$\mathbb R$ is a limit of points of the form $(m+\alpha_j/2)/q_j\in H_j$,
since $q_j\ge2^j/\operatorname{diam}A$ by \eqref{SN:eq:q-lower}. Hence
$\bigcup_jH_j$ is dense and $E$ is nowhere dense.
\\
\noindent\textbf{Step 3: Bound the length and number of components of $D_E(A)$ .} Let $t\in D_E(A)$ and
pick $x$ with $x+tA\subseteq E$. For all $1\le i\le n_j$, the point
$x+t(x_j+\sigma_jb_{j,i})$ lies outside $H_j$, which means
\[q_jx+q_jtx_j+ q_jt\sigma_jb_{j,i} (\bmod 1)\notin (0,\alpha_j).\]
This is equivalent to saying that
$q_j t\sigma_j\,b_{j,i}\bmod 1$ avoid the arc
$\bigl(-q_j(x+tx_j),\,\alpha_j-q_j(x+tx_j)\bigr)\bmod1$ of length
$\alpha_j$. By Lemma~\ref{lem:missing-arc},
$\sigma_jt\in B(P_j,q_j)$, and therefore
\[
 D_E(A)\cap[-j,j]\subseteq
 B_j:=\bigl\{t\in[-j,j]:\sigma_jt\in B(P_j,q_j)\bigr\}.
\]
Apply
Lemma~\ref{lem:missing-arc} with $\nu$ Lebesgue measure on
$[-j,j]$, for which $|\widehat\nu(\xi)|\le1/|\xi|$. By
\eqref{SN:eq:q-lower},
\[
 \nu(B_j)\le C_j\Bigl[2je^{-k_jn_j}+\frac{W_j^{n_j}}{q_jb_{j,n_j}}\Bigr]
 \le C_j(2j+1)e^{-k_jn_j}=:L_j.
\]
Since $\sum_ib_{j,i}/b_{j,n_j}\le\frac43e^{S_j}$ by Lemma~\ref{lemma:tail-domination}, the component bound of
Lemma~\ref{lem:missing-arc} and the bound $q_j\le2W_j^{n_j}e^{k_jn_j}/b_{j,n_j}$ following \eqref{SN:eq:q-lower} tell us that the number of connected components of $B_j$ is at most
\[
C_j\Bigl(n_j+jq_jb_{j,n_j}\cdot\tfrac43e^{S_j}\Bigr)\leq
C_j\Bigl(n_j+j\cdot(2W_j^{n_j}e^{k_jn_j})\cdot\tfrac43e^{S_j}\Bigr)
 \le e^{S_j+(D_j-k_j)n_j}=:N_j,
\]
using $n_j\le e^{n_j}$ and the definition of $D_j$.
\\ \\
\noindent\textbf{Step 4: Construct the open cover for the gauge Hausdorff measure bound.} Let $d_j=L_j/N_j$. Cutting each component $J$
of $B_j$ into $\lceil|J|/d_j\rceil\le1+|J|/d_j$ intervals of length at
most $d_j$ covers $B_j$ by at most $N_j+L_j/d_j=2N_j$ intervals. Moreover
$d_j\le L_j\le1$ and
\[
 \log\frac1{d_j}=S_j+(D_j-k_j)n_j+k_jn_j-\log\bigl(C_j(2j+1)\bigr)
 \le S_j+D_jn_j.
\]
By \eqref{SN:eq:overhead-bound}, monotonicity of $\Psi_h$, and
\eqref{SN:eq:block-choice}, the sum $\sum_U h(\operatorname{diam}U)$ over the intervals $U$ of this cover is at most
\[
2N_jh(d_j)\le2N_jd_j\,e^{\Psi_h(\log(1/d_j))}
 \le2L_j\,e^{k_jn_j/2}
 =2C_j(2j+1)e^{-k_jn_j/2}<2^{-j}.
\]
Fix an integer $T\ge1$. For every $j\ge T$, Step~3 gives $D_E(A)\cap[-T,T]\subseteq B_j$, and we just constructed a cover of $B_j$ by intervals of length at most $d_j\le2^{-j}$ with total $h$-cost below $2^{-j}$. 
Hence $\mathcal{H}^h(D_E(A)\cap[-T,T])=0$, and countable
subadditivity over $T$ gives $\mathcal{H}^h(D_E(A))=0$.
\end{proof}
We conclude this section with the proof of Proposition~\ref{SN:prop:pairs}.

\begin{proof}[Proof of Proposition~\ref{SN:prop:pairs}]
Pick two distinct points $a,a'\in A$. If $x+tA\subseteq K\subseteq[0,1]$,
then $t(a-a')=(x+ta)-(x+ta')\in[-1,1]$, so $|t|\le|a-a'|^{-1}$.
Fixing also $a_0\in A$ and using $x+ta_0\in[0,1]$, the pair $(x,t)$
is confined to a rectangle of dimensions $X\times2|a-a'|^{-1}$, where
$X=1+2|a_0|/|a-a'|$.

Fix $\eta>0$. Since $\mathcal{H}^h(D_K(A))=0$, cover $D_K(A)$ by
intervals $I_1,I_2,\ldots$ of lengths $r_i<\min\{1,\eta/2\}$ with
$\sum_ih(r_i)<\eta$. The slice of $\mathcal{E}_K(A)$ over each $I_i$
is contained in a rectangle of dimensions $r_i\times X$, which can be
covered by at most $C_X/r_i$ squares of side $r_i$ and diameter $2r_i<\eta$.
Their total $g_h$-cost is
\[
 \sum_i\frac{C_X}{r_i}\cdot g_h(2r_i)
 =\sum_i\frac{C_X}{r_i}\cdot 2r_ih(r_i)
 =2C_X\sum_ih(r_i)<2C_X\eta.
\]
Since $\eta>0$ was arbitrary, $\mathcal{H}^{g_h}(\mathcal{E}_K(A))=0$.
That $g_h$ is a gauge follows from $u$ and $h(u/2)$ both being
nondecreasing and tending to zero at the origin.
\end{proof}
\subsection{Examples and corollaries}
The hypothesis of Theorem~\ref{SN:thm:block-entropy} is easy to check for explicit sequences.

\begin{corollary}\label{SN:cor:exact-thresholds}
Let $A\subseteq\mathbb R$ be bounded and contain affine copies of
arbitrarily long initial segments of one of the sequences $(b_n)$
below, with affine maps allowed to vary with the length. Let
$(h_m)_{m\ge1}$ be a countable family of gauges,
each satisfying the corresponding condition on their overhead functions:
\begin{enumerate}
\item $b_n=1/n!$, \;\;\;$\Psi_{h_m}(S)=o(S/\log S)$;
\item $b_n=e^{-n^q}$, $q>1$, \;\;\; $\Psi_{h_m}(S)=o(S^{1/q})$;
\item $b_n=e^{-e^{n^\beta}}$, $\beta>0$ , \;\;\;
      $\Psi_{h_m}(S)=o\bigl((\log S)^{1/\beta}\bigr)$.
\end{enumerate}
Then for every $\varepsilon>0$ there is a closed $1$-periodic
nowhere-dense $E$ with $\mathcal L^1(E\cap I)>1-\varepsilon$ in every unit
interval and $\mathcal{H}^{h_m}(D_E(A))=0$ for every $m\ge1$;
$K=E\cap[0,1]$ gives the compact conclusion.
\end{corollary}

\begin{proof}
Fix $m$ and $\rho\in(0,1)$, and choose $N_\rho$ with $b_{n+1}/b_n\le\rho$ for all
$n\ge N_\rho$. For every $n$ the set $A$ contains an affine copy
$x+t\{b_{N_\rho},\ldots,b_{N_\rho+n-1}\}$. Writing $\sigma=\operatorname{sgn}t$
and $\tilde b_i=|t|\,b_{N_\rho+i-1}$, this is a set $x+\sigma\{\tilde b_1,\ldots,\tilde b_n\}$
with $\tilde b_{i+1}\le\rho\,\tilde b_i$ and $\log(\tilde b_1/\tilde b_n)=S_\rho(n)$, where
$S_\rho(n)=\log(b_{N_\rho}/b_{N_\rho+n-1})$. Direct computation gives
\begin{align*}
 S_\rho(n)&=n\log n+O_\rho(n)
 &&\text{(case 1)},\\
 S_\rho(n)&=n^q+O_{\rho,q}(n^{q-1})
 &&\text{(case 2)},\\
 \bigl(\log S_\rho(n)\bigr)^{1/\beta}&=n+O_{\rho,\beta}(1)
 &&\text{(case 3)}.
\end{align*}
In each case the stated condition gives $\Psi_{h_m}(S_\rho(n)/\rho)=o(n)$, so
$\Psi_{h_m}(S_\rho(n)/\rho)\le\rho n$ for all large $n$. Hence
\eqref{SN:eq:block-criterion} holds, and
Theorem~\ref{SN:thm:block-entropy} applies.
\end{proof}

Explicit gauges satisfying each condition include the following.
For $b_n=1/n!$: take $h(u)=u\exp\!\bigl(S/(\log S\log\log S)\bigr)$
with $S=\log(1/u)$. For $b_n=e^{-n^q}$: any gauge with overhead
$\Psi_h(S)=S^{1/q}/\log S$. For $b_n=e^{-e^{n^\beta}}$: the families
\[
 \widetilde h_\gamma(u)
 =u\exp\!\bigl((\log\log(1/u))^\gamma\bigr),
 \quad0<\gamma<1/\beta,
 \qquad
 \widetilde\ell_p(u)=u(\log\log(1/u))^p,
 \quad p>0.
\]

\begin{remark}
In each case, one set $E$ can be chosen to work simultaneously for
the entire indicated range of parameters. For $b_n=e^{-n^q}$, choose
$\theta_j\uparrow1/q$ and apply Corollary~\ref{SN:cor:exact-thresholds}
to $(h_{\theta_j})_{j\ge1}$. Every $h_\theta$ with $\theta<1/q$ is
dominated near zero by some $h_{\theta_j}$, so
$\mathcal{H}^{h_\theta}(D_E(A))=0$ for all such $\theta$. The same
argument applies to the doubly exponential case by taking
$\gamma_j\uparrow1/\beta$ and positive integers for $\widetilde\ell_p$.
\end{remark}

\begin{corollary}\label{cor:joint-conclusions}
Let $(A_\ell)_{\ell\ge1}$ be a countable family of bounded infinite
sets, let $(\omega_m)_{m\ge1}$ be decay moduli, and for each $\ell$
let $(h_m^{(\ell)})_{m\ge1}$ be a countable family of gauges
satisfying \eqref{SN:eq:block-criterion} for $A=A_\ell$. Then for every $\varepsilon>0$
there is a single closed $1$-periodic nowhere-dense $E$ with
$\mathcal L^1(E\cap I)>1-\varepsilon$ in every unit interval such that
$\nu(D_E(A_\ell))=0$ for every $\ell$ and every finite positive
measure $\nu$ with decay $\omega_m$ for some $m$, and
$\mathcal{H}^{h_m^{(\ell)}}(D_E(A_\ell))=0$ for every $\ell$ and $m$.
\end{corollary}

\begin{proof}
Apply Theorem~\ref{thm:prescribed-moduli} to each $A_\ell$ with
budget $\varepsilon_\ell^F$, and Theorem~\ref{SN:thm:block-entropy}
to each $A_\ell$ and gauge family $(h_m^{(\ell)})$ with budget
$\varepsilon_\ell^G$, choosing all budgets so that
$\sum_\ell(\varepsilon_\ell^F+\varepsilon_\ell^G)<\varepsilon$.
Let $E$ be the intersection of all resulting sets. Closedness and
periodicity are preserved; for each $\ell$, $E$ is contained in the
set produced for $A_\ell$, so $D_E(A_\ell)$ is contained in that
set's exceptional parameter set, and all nullity conclusions hold.
The measure satisfies $\mathcal L^1(E\cap I)\ge1-\sum_\ell(\varepsilon_\ell^F+\varepsilon_\ell^G)>1-\varepsilon$.
\end{proof}


\end{document}